\documentclass[preprint,times,12pt]{article}

\usepackage{amstext}
\usepackage{latexsym}
\usepackage{amssymb}
\usepackage{amsmath}
\usepackage{graphicx}
\usepackage{ifpdf}
\usepackage{color}

\allowdisplaybreaks
\DeclareGraphicsExtensions{.pdf,.png,.jpg,.jpeg,.mps,.eps}

\usepackage[colorlinks=false,linktocpage=true]{hyperref}
\usepackage{hyperref}
\hypersetup{colorlinks=false}

\def\R{\mathbb R}

\def\C{\mathbb C}
\def\N{\mathbb N}

\def\G{\mathsf{G}}

\def\int{{\rm int\;}}

\def\tr{{\rm tr}}
\def\eps{\varepsilon}
\def\PSL{{\rm PSL}}

\def\X{{\mathfrak X}}
\def\H{\mathbb H}

\def\X{\mathcal X}

\newcommand{\SL}{{\rm SL}}

\newtheorem{theorem}{Theorem}[section]

\newtheorem{definition}[theorem]{Definition}

\newtheorem{proposition}[theorem]{Proposition}

\newtheorem{remark}[theorem]{Remark}

\newtheorem{lemma}[theorem]{Lemma}

\usepackage{hyperref}

\begin{document}
		\title{Expansiveness for the geodesic and horocycle flows on compact Riemann surfaces
		of constant negative curvature}
	\author{{\sc Huynh Minh Hien} \\[1ex] 
		Department of Mathematics,\\
		Quy Nhon University,\\
		 170 An Duong Vuong, Quy Nhon, Vietnam;\\[1ex] 
		e-mail: huynhminhhien@qnu.edu.vn \date{} 
	}
	\maketitle 
	\begin{abstract}We study expansive properties for the geodesic and horocycle flows on  Riemann surfaces 
		of constant negative curvature. It is well-known that the geodesic flow is expansive in the sense of Bowen-Walters
		and the horocycle flow is positive and negative strong separating in the sense of Gura. 
		In this paper, we give a new proof for the expansiveness in the sense of Bowen-Walters for the geodesic flow
		and show that the horocycle flow is positive and negative strong kinematic expansive in the sense of Artigue
		as well as expansive in the sense of Katok-Hasselblatt but not expansive in the sense of Bowen-Walters. We also point out that the geodesic flow is neither positive nor negative separating.	
	\end{abstract}

	
	\maketitle
	
	\section{Introduction}
	The study of expansive flows started in 1972 with the works of Bowen/Walters \cite{bw} and Flinn \cite{flinn}. In \cite{bw}, the authors  generalized the definition of expansive homeomorphisms to introduce a reasonable definition of expansiveness for flows  that is called `expansive in the sense of Bowen and Walters' (or shortly BW\,--\,expansive). Since then, 
	there have been different varieties of expansive flows introduced. 
	In 1984, Komuro \cite{komuro} gave the notions `C\,--\,expansive' (as the same to `BW\,--\,expansive'), `K\,--\,expansive' (as the same to `expansive' in the same of Flinn) and  `K$^*$\,--\,expansive' to investigate geometric Lorenz attractors.
	In general, `K\,--\,expansive' is weaker than `BW\,--\,expansive' but stronger than `K$^*$\,--\,expansive'.
	In the case of fixed-point-free flows on compact metric spaces, the three notions  are equivalent (see \cite{artigue2,oka}).

	A different and very interesting kind of expansiveness  called `separating' was discovered by Gura \cite{gura} in 1984. The author showed that the horocycle flow on a compact surface with negative curvature is positive and
	negative separating.  His definition in \cite{gura} requires to separate every pair of points in
	different orbits. The author also proved a remarkable result: every global time change of such
	flow is positive and negative separating (such property is also called positive and negative strong separating). In 1995, Katok/Hasselblatt \cite{kh} gave another kind of expansiveness (also called KH\,--\,expansiveness) which is weaker than BW\,--\,expansiveness but implies separation. 
	It then was showed by Artigue \cite{artigue2} that	a flow is KH\,--\,expansive if and only if it is separating and the set of its fixed points is open.

	Recently, in 2016, Artigue \cite{artigue} used the term   `geometric expansive' as {\em K}\,--\,expansive and introduced the term `kinematic expansive' which
	is a stronger property than separation and weaker than BW\,--\,expansiveness.
	The author also considered the forms of `strong kinematic expansive', `geometric separating', `strong separating' and `separating' flows.  Examples are given to analyze the relationships among the above definitions. Some interesting properties are proved in different contexts: surfaces, suspension flows and compact metric spaces.

	Regarding properties of the geodesic flow, in 1967, Anosov \cite{anosov} showed that the geodesic flow on compact Riemannian manifolds with negative curvature 
	is hyperbolic. In 1972, it was  proved by Bowen \cite{bowen-pe} that  hyperbolic flows are BW\,--\,expansive
	and consequently the geodesic flow on compact Riemannian manifolds of negative curvature is BW\,--\,expansive.
	As mentioned above, in 1984, Gura \cite{gura} showed that the horocycle on compact surface of negative curvature is positive and negative separating.

The paper is organized as follows. In the next section, we introduce the necessary background 
material which is well-known in principle \cite{bedkeanser,einsward,ratcliff}.  
 Section \ref{sec3} is devoted to consider expansive properties mentioned above for the geodesic and horocycle flows on compact factors of the hyperbolic plane. A new detailed proof for BW\,--\,expansiveness of the geodesic flow (Theorem \ref{geoex}) via a  property of the injectivity radius is given. The horocycle flow on this model is not only positive/negative strong separating but also positive/negative strong kinematic expansive  (Theorem \ref{pne}) as well as KH\,--\,expansive (Theorem \ref{khthm}). In the end, we point out that the horocycle flow is not BW\,--\,expansive and the geodesic flow is neither positive  nor negative separating; see Remark \ref{rm}. 
	\section{Preliminaries}
	
	We consider the geodesic and horocycle flows on compact Riemann surfaces of constant negative curvature. 
	It is well-known that any compact orientable surface 
	with constant negative curvature is isometric to a factor $\Gamma\backslash \H^2=\{\Gamma z, z\in\H^2 \}$, 
	where $\H^2=\{z=x+iy\in \C:\, y>0\}$ is the hyperbolic plane endowed  
	with the hyperbolic metric $ds^2=\frac{dx^2+dy^2}{y^2}$ and
	$\Gamma $ is a Fuchsian group that is  discrete subgroup of the group $\PSL(2,\R)=\SL(2,\R)/\{\pm E_2\}$; 
	where $\SL(2,\R)$ is the group of all real $2\times 2$ matrices with unity determinant, 
	and $E_2$ denotes the unit matrix. In the hyperbolic plane model,
	geodesics are vertical lines and semi-circles centered on the real axis.
	 The group $\PSL(2,\R)$  acts transitively on $\H^2$ by 
	M\"obius transformations
	$z\mapsto \frac{az+b}{cz+d}$.
	If the action is free of fixed points, then the factor $\Gamma\backslash\H^2$  has a Riemann surface structure
	that is a closed Riemann surface of genus at least $2$ 
	and has the hyperbolic plane $\H^2$ as the universal covering.
	The unit tangent bundle $T^1\H^2$ is isometric to the group $\PSL(2,\R)$ 
	and as a consequence,  the unit tangent bundle $T^1(\Gamma\backslash\H ^2)$
	is isometric to the quotient space 
	$\Gamma\backslash \PSL(2,\R)=\{\Gamma g,g\in\PSL(2,\R)\}$, 
	which is the system of right co-sets of $\Gamma$ in $\PSL(2,\R)$, by an isometry
	$\Xi$. 	Since ${\rm PSL}(2, \R)$ 
	is connected, also $\Gamma\backslash {\rm PSL}(2, \R)$ is connected. 
	Furthermore, $X=\Gamma\backslash {\rm PSL}(2, \R)$
	is a three-dimensional real analytic manifold.
		
	The geodesic flow on $T^1\H^2$ can be described as the flow
	$\varphi_t^\G(g)=ga_t$ on $\G:=\PSL(2,\R)$, where $a_t\in\G$ denotes the equivalence class
	 obtained from the matrix $A_t=\scriptsize\big(\begin{array}{cc}
	e^{t/2} & 0\\ 0 & e^{-t/2}
	\end{array}\big)$, and whence the geodesic flow $(\varphi_t^\X)_{t\in\R}$ on $\X=T^1(\Gamma\backslash\H^2)$
	can be described as the `quotient flow' 
	\begin{equation*}\varphi^X_t(\Gamma g)=\Gamma g a_t
	\end{equation*}
	on $X=\Gamma\backslash\PSL(2,\R)$  by the conjugate relation
	\begin{equation}\label{ce}
	\varphi_t^\X=\Xi^{-1}\circ\varphi_t^X\circ\Xi.
	\end{equation} 
	
   A horocycle is a (euclidean) circle tangent to real axis or a horizontal line. The 
   stable and unstable horocycle flows on $T^1\H^2$ can be described as the flows: $\theta_t^\G(g)=gb_t,\, \eta_t^\G(g)=gc_t$ on $\G$; 
   where
   $b_t,c_t\in\PSL(2,\R)$ denote the equivalence classes obtained from the matrices $B_t=\scriptsize\big(\begin{array}{cc}
   	1 & t\\ 0 & 1
   \end{array}\big),C_t=\scriptsize\big(\begin{array}{cc}1 &0\\ t&1 \end{array} \big)\in\SL(2,\R)$. 
   Therefore the  stable and unstable horocycle flows $(\theta_t^\X)_{t\in\R}$, $(\eta_t^\X)_{t\in\R}$  on $\X=T^1(\Gamma\backslash\H^2)$ can be equivalently described as the flows
	\begin{equation*} 
	\theta^X_t(\Gamma g)=\Gamma g b_t,
	\quad \eta^X_t(\Gamma g)=\Gamma g c_t
	\end{equation*} 
	on $X=\Gamma\backslash\PSL(2,\R)$ 
	by the conjugate relations 
	\begin{equation}\label{cr}  
	\theta_t^\X=\Xi^{-1}\circ\theta_t^X\circ\Xi,\quad \eta_t^\X=\Xi^{-1}\circ\eta_t^X\circ\Xi \quad \mbox{for all}\quad t\in\R.
	\end{equation}
	
	\smallskip 
	There are some more advantages to work on $X=\Gamma\backslash\PSL(2,\R)$
	rather than on $\X=T^1(\Gamma\backslash\H^2)$. For example, one can calculate explicitly the stable and unstable manifolds 
	at a point $x$ to be
	\begin{equation*}
	W^s_X(x)=\{\theta^X_t(x),t\in\R\}
	\quad \mbox{and}\quad W^u_X(x)=\{\eta^X_t(x), t\in\R\}.
	\end{equation*} 
	The flow $(\varphi^X_t)_{t\in\R}$
	is hyperbolic, that is, 
	for every $x\in X$ there exists an orthogonal and $(\varphi_t^X)_{t\in\R}$-stable splitting of the tangent space
	$T_xX$
	\[T_x X= E^0(x)\oplus E^s(x)\oplus E^u(x)\]
	such that the differential of the flow $(\varphi_t^X)_{t\in \R}$ is uniformly expanding on $E^u(x)$,
	 uniformly contracting on $E^s(x)$ and isometric on $E^0(x)=\langle \frac{d}{dt}\varphi_t^X(x)|_{t=0}\rangle$. One can choose  
	\begin{eqnarray*}
		E^s(x)  =  \Big\langle\frac{d}{dt}\,\theta^X_t(x)\Big|_{t=0}\Big\rangle
		\quad\mbox{and}\quad 
		E^u(x)  =  \Big\langle\frac{d}{dt}\,\eta^X_t(x)\Big|_{t=0}\Big\rangle.
	\end{eqnarray*}

The horocycle flows $(\theta^\X)_{t\in\R}$ and $(\eta^\X)_{t\in\R}$ are ergodic \cite{markus}. 
If the space $\Gamma\backslash\H^2$ has a finite volume,  each orbit is either periodic or
dense.
 In the case that the space $\Gamma\backslash\H^2$  is compact,
 there are no periodic orbits for the horocycle flows.

	General references for this section are \cite{bedkeanser,einsward}, 
	and these works may be consulted for the proofs to all results which are stated above.
	In what follows,
	 we will drop the superscript $X$ from $(\varphi^X_t)_{t\in\R},(\theta^X_t)_{t\in\R},(\eta^X_t)_{t\in\R}$ to simplify 	notation. 
	We consider the stable horocycle flow only and use the term `horocycle flow' for it. In the whole present paper, we always assume the action of $\Gamma$ on $\H^2$ to be free (of fixed points) and the factor   $\Gamma\backslash\H^2$  to be compact. Note that $\Gamma\backslash\H^2$ is compact if and only if $\Gamma\backslash\PSL(2,\R)$ is compact.

	In the rest of this section  we collect some notions and useful technical results. 
	
		\begin{lemma}\label{at}
		There is a natural Riemannian metric on $\G=\PSL(2,\R)$
		such that the induced metric function $d_\G$ is left-invariant under $\G$ and 
		\[d_\G(a_t,e)=\frac{1}{\sqrt 2}|t|, \quad d_\G(b_t,e)\leq |t|,\quad d_\G(c_t,e)<|t|\quad\mbox{for all}\quad t\in\R \]
		where $e=\pi(E_2)$ is the unity of $\G$.
	\end{lemma}	

	We define a metric function $d_{X}$ on $X=\Gamma\backslash\PSL(2,\R)$ by 
	\[ d_{X}(x_1, x_2)
	=\inf_{\gamma_1, \gamma_2\in\Gamma} d_{\G}(\gamma_1 g_1, \gamma_2 g_2)
	=\inf_{\gamma\in\Gamma} d_{\G}(g_1, \gamma g_2), \]   
	where $x_1=\Gamma g_1 $, $x_2=\Gamma g_2$.
	In fact, if $X$ is compact, one can prove that the infimum is a minimum:
	\begin{equation*}
	d_{X}(x_1, x_2)=\min_{\gamma\in\Gamma} d_{\G}(g_1, \gamma g_2).
	\end{equation*} 
		It is possible to derive 
	a uniform lower bound on $d_{\G}(g, \gamma g)$ 
	for $g\in \PSL(2,\R)$ and $\gamma\in\Gamma\setminus\{e\}$.
	
	\begin{lemma}\label{sigma_0}
		If the space $X=\Gamma\backslash\PSL(2,\R)$ is compact, then there exists $\sigma_0>0$ such that
		\[d_\G(\gamma g, g)>\sigma_0\quad\mbox{for all}\quad \gamma\in\Gamma\setminus\{e\}. \]
	\end{lemma}
	The number $\sigma_0$ is called an injectivity radius. See \cite[Lemma 1, p. 237]{ratcliff} for a similar result. 
	
	For $g=\pi(G)\in\PSL(2,\R), G=\big(\scriptsize\begin{array}{cc}a&b\\c&d\end{array} \big)$,
	the trace of $g$ is defined by \[\tr(g)=|a+d|.\]
	If the action of $\Gamma$ on $\H^2$ is free and the factor $\Gamma\backslash\H^2$ is compact
	then all elements $g\in\Gamma\setminus\{e\}$ are hyperbolic \cite[Theorem 6.6.6]{ratcliff}, i.e. $\tr(g)>2$. 
	Furthermore, one gets a stronger result:
	\begin{lemma}\label{hy}
		If the factor $\Gamma\backslash\H^2$ is compact, then there exists $\eps_*>0$ such that
		\[\tr(g)\geq 2+\eps_*\quad \mbox{for all}\quad g\in \Gamma\setminus \{e\}. \]
	\end{lemma}	
	
	 Here are some more auxiliary results.
	
	\begin{lemma}\label{konvexa}
		(a) For every $\delta>0$ there is $\rho>0$ with the following property. 
		If $G=\Big(\begin{array}{cc}g_{11}&g_{12}\\
		g_{21}&g_{22}
		\end{array}\Big)\in\SL(2,\R)$ satisfies
		$|g_{11}-1|+|g_{12}|+|g_{21}|+|g_{22}-1|<\rho$ then $d_{\G}(g, e)<\delta$ for $g=\pi(G)$, 
		where $\pi: {\rm SL}(2, \R)\to {\rm PSL}(2, \R)$ is the natural projection.
		
		\smallskip
		\noindent
		(b) For every $\eps>0$ there is $\delta>0$ with the following property. 
		If $g\in\G$ satisfying $d_\G(g,e)<\delta$ then there are
		\[ G=\Bigg(\begin{array}{cc}g_{11}&g_{12}\\
		g_{21}&g_{22}
		\end{array}\Bigg)\]
		such that $g=\pi(G)$ and
		\[|g_{11}-1|+|g_{12}|+|g_{21}|+|g_{22}-1|<\eps.\]
	\end{lemma}

{\bf Proof\,:} (a) See \cite[Lemma 2.17 (a)]{HK} for a proof.

(b) Indeed, suppose on contrary that  
\begin{equation}\label{hodett} 
|g^j_{11}-1|+|g^j_{12}|+|g^j_{21}|+|g^j_{22}-1|\ge\eps_0
\end{equation} 
for some sequence $d_{\G}(g^j, e)\to 0$ and all $G^j\in {\rm SL}(2, \R)$ 
such that $g^j=\pi(G^j)$. For $j\in\N$ take any $G^j\in\SL(2,\R)$ so that $g^j=\pi(G^j)$. 
From (a) we deduce that $|g^j_{12}|+|g^j_{21}|\to 0$, $|g^j_{11}|\to 1$, 
$|g^j_{22}|\to 1$, and $g^j_{11} g^j_{22}\to 1$. Thus, along a subsequence 
which is not renamed, either $g^j_{11}\to 1$, $g^j_{22}\to 1$ 
or $g^j_{11}\to -1$, $g^j_{22}\to -1$. The first case is impossible in view of (\ref{hodett}). 
In the second case we consider $\tilde{G}^j=-G^j$ which also has $g^j=\pi(\tilde{G}^j)$. 
But then (\ref{hodett}) implies 
\[ |g^j_{11}+1|+|g^j_{12}|+|g^j_{21}|+|g^j_{22}+1|\ge\eps_0, \] 
and once more this is impossible. 

\begin{definition}\rm
	Let $\phi:\R\times M\to M$ be a flow. 
	
	\begin{enumerate}
		\item[(a)] A point $x\in M$ is called a {\em fixed point} (or {\em singular point})
		if \[\phi_t(x)=x\quad\mbox{for all}\quad t\in\R. \]
		\item[(b)] A point $x\in M$ is called a {\em periodic point} if there is $T>0$ such that
		\[\phi_T(x)=x. \]
	\end{enumerate} 
\end{definition}

\begin{proposition}\label{fixedpoint} Assume that $X=\Gamma\backslash\PSL(2,\R)$ is compact.
	Then the flow $(\theta_t)_{t\in\R}$ does not have a periodic point. 
	In particular, it has no fixed points. 
\end{proposition}
\noindent
{\bf Proof\,:} Suppose in contrary that $x=\Gamma g$ is a periodic point of $(\theta_t)_{t\in\R}$,
i.e. $\theta_T(x)=x$ for some $T>0$. Then $g^{-1}\gamma g=b_T$ for some $\gamma\in\Gamma$ implies $\tr(\gamma)=\tr(b_T)=2$. It follows from  Lemma \ref{hy} that  $\gamma=e$. Therefore $g=gb_T$ yields $T=0$ which is a contradiction. The latter assertion is obvious.
{\hfill $\Box$}

\begin{definition}[\cite{artigue}]\rm
	Two continuous flows $\phi:\R\times X\to X$ and $\psi: \R\times Y\to Y$ are said to be 
	{\em equivalent} if there is a homeomorphism $h:X\to Y$ such that $\phi_t=h^{-1}\psi_th$ for all $t\in\R$. 
\end{definition}
Via \eqref{ce} and \eqref{cr}, the flows $(\varphi_t)_{t\in\R}$ and $(\varphi^\X_t)_{t\in\R}$ are equivalent, 
and so are $(\theta^\X_t)_{t\in\R}$ and $(\theta_t)_{t\in\R}$.  
It is easy to see that all the expansive properties introduced in the next section are invariant under equivalence.  

\section{Expansive properties}\label{sec3}
In this section we study BW\,--\,expansive,  kinematic expansive, separating, and KH\,--\,expansive properties for 
the geodesic flow $(\varphi_t^\X)_{t\in\R}$ and the horocycle flow $(\theta^\X_t)_{t\in\R}$ on $\X=T^1(\Gamma\backslash\H^2)$. 
We reprove that the geodesic flow is BW\,--\,expansive. The horocycle flow is positive/negative strong kinematic expansive
as well as KH\,--\,expansive but not BW\,--\,expansive. 
\subsection{BW\,--\,expansiveness}
This subsection provides a new detailed proof of the expansiveness in the sense of Bowen-Walters
for the geodesic flow $(\varphi_t^\X)_{t\in\R}$ owing to a characteristic property of the injectivity radius. 

\begin{definition}[\cite{bw}, BW\,--\,expansive]\rm Let $(M,d)$ be a compact metric space. 
	A continuous flow $\phi:\R\times M\longrightarrow M$ is called {\em BW\,--\,expansive} if for each $\eps>0$, 
	there exists $\delta>0$ with the following property. 
	If $s:\R\rightarrow \R$ is a continuous function with $s(0)=0$ and 
	\[d(\phi_t(x),\phi_{s(t)}(y))<\delta \quad \mbox{for all}\quad t\in \R\]
	then $y=\phi_\tau(x)$ for some $\tau\in (-\eps,\eps)$.
\end{definition}
It was showed in \cite[Theorem 3]{bw} that `continuous function' in the above definition can be replaced by
`increasing homeomorphism' in the case of fixed-point-free flows. 
Then this definition becomes the one in \cite{flinn}, `K\,--\,expansive' \cite{komuro}, and
`geometric expansive' \cite{artigue}.  
\begin{theorem}\label{geoex}
	The geodesic flow $(\varphi^\X_t)_{t\in\R}$ is BW\,--\,expansive.
\end{theorem}
\noindent{\bf Proof\,:} Since BW\,--\,expansiveness is an invariant under equivalence 
\footnote{It is showed in \cite[Corollary 4]{bw} that BW\,--\,expansiveness is an invariant under conjugacy that is weaker than equivalence. Recall that the flows $(\phi_t)_{t\in\R}$ on $X$  and $(\psi_t)_{t\in\R}$ on $Y$ are said to be conjugate if
there is a homeomorphism from $X$ to $Y$ mapping the orbits of $(\phi_t)_{t\in\R}$
onto orbits of $(\psi_t)_{t\in\R}$.}, it follows from \eqref{ce} that it suffices to show that
the flow $(\varphi_t)_{t\in\R}$ is BW\,--\,expansive. Let $\eps>0$ be given, 
$\eps_0=e^{\eps/2}-e^{-\eps/2}>0$ and set $\delta=\delta(\eps_0)<\sigma_0/4$ as in Lemma \ref{konvexa}\,(b);
here $\sigma_0$ is from Lemma \ref{sigma_0}. Let $x,y\in X$ and $s:\R\to\R$ be continuous with $s(0)=0$ such that
\begin{equation*}
d_X(\varphi_{s(t)}(y),\varphi_t(x))<\delta\quad \mbox{for all}\quad t\in\R.
\end{equation*}
Write $x=\Gamma g, y=\Gamma h$ for $g,h\in\G$. 
For every $t\in \R$, there is $\gamma(t)\in\Gamma$ so that
\begin{eqnarray}\label{gammat}
d_X(\varphi_{s(t)}(y),\varphi_t(x))
=d_X(\Gamma ha_{s(t)}, \Gamma ga_t)
=d_\G(ha_{s(t)},\gamma(t) ga_t)<\delta. 
\end{eqnarray}
We claim that  $\gamma(t)=\gamma(0)=:\gamma$ for all $t\in\R$. 
For any $t_1,t_2\in\R$, we have
\begin{eqnarray*}
	\lefteqn{d_\G(\gamma(t_2)^{-1}\gamma(t_1)ga_{t_1},ga_{t_1})}
	\\
	&=& d_\G(\gamma(t_1)ga_{t_1},\gamma(t_2)ga_{t_1})
	\\
	&\leq& d_\G(\gamma(t_1)ga_{t_1},ha_{s(t_1)})
	+d_\G(ha_{s(t_1)},ha_{s(t_2)})
	+d_\G(ha_{s(t_2)},\gamma(t_2)ga_{t_2})
	\\
	&&+\ d_\G(\gamma(t_2)ga_{t_2},\gamma(t_2)ga_{t_1})
	\\
	&=& d_\G(\gamma(t_1)ga_{t_1},ha_{s(t_1)})
	+d_\G(a_{s(t_1)},a_{s(t_2)})
	+d_\G(ha_{s(t_2)},\gamma(t_2)ga_{t_2})
		+d_\G(a_{t_2},a_{t_1} )
	\\
	&\leq& 2\delta+\frac{1}{\sqrt 2}|s(t_1)-s(t_2)|+\frac{1}{\sqrt 2}|t_1-t_2|,
\end{eqnarray*}
due to Lemma \ref{at}.
For given $L>0$, we verify that $\gamma(t)=\gamma(0)$ for all $t\in [-L,L]$. Indeed, since $s:[-L,L]\to\R$ is uniformly continuous, there is $0<\rho=\rho(L,\delta)<\delta$
such that if $t_1,t_2\in [-L,L]$ and $|t_1-t_2|<\rho$ then 
$|s(t_1)-s(t_2)|<\delta$.
For $t_1,t_2\in [0,\rho/2]$, then $|t_1-t_2|<\rho$ implies $|s(t_1)-s(t_2)|<\delta$. 
This yields 
\[d_\G(\gamma(t_2)^{-1}\gamma(t_1)c_1(t_1),c_1(t_1))<4\delta <\sigma_0.\]
From the property of $\sigma_0$ in Lemma \ref{sigma_0}, it follows that 
$\gamma(t_2)=\gamma(t_1)$ for $|t_1-t_2|<\rho$. 
Here if we specialize this to $t_1=0$ and $t_2\in [0,\rho/2]$,
then $\gamma(t_2)=\gamma(0)$ for all $t_2\in [0,\rho/2]$. 
 Then we repeat the argument for
$t_1=\rho/2$ and $t_2\in [\rho/2, \rho]$, we deduce that $\gamma(t)=\gamma(0)$ for all $t\in [0,\rho]$, which upon further iteration leads to $\gamma(t)=\gamma(0)$ for all $t\in [0,L]$
and similarly $\gamma(t)=\gamma(0)$ for all $t\in [-L,0]$. Therefore,  
\begin{equation}\label{st}
d_X(\varphi_{s(t)}(y),\varphi_t(x))=d_\G(a_{-t}g^{-1}\gamma ha_{s(t)},e)<\delta \quad\mbox{for all}\quad t\in\R.
\end{equation}
Write $g^{-1}\gamma h=\pi(K)$ for $K=\big(\scriptsize\begin{array}{cc} a&b\\c&d\end{array}\big)\in\SL(2,\R)$.
Thus
\[A_{-t}KA_{s(t)}=\Big(\begin{array}{cc} ae^{\frac{s(t)-t}2}&be^{-\frac{s(t)+t}2}\\c e^{\frac{s(t)+t}2}&d e^{\frac{t-s(t)}2}\end{array}\Big)\]
together with \eqref{st} implies 
\begin{equation}\label{abcdeq}
\big||a|e^{\frac{s(t)-t}2}-1\big|+|b|e^{-\frac{s(t)+t}2}
+|c| e^{\frac{s(t)+t}2}+
\big||d| e^{\frac{t-s(t)}2}-1\big|<\eps_0\quad \mbox{for all}\quad t\in\R,
\end{equation} 
using Lemma \ref{konvexa}\,(b).
Then there is $M>0$ such that $|s(t)-t|\leq M$ for all $t\in\R$ and hence
$s(t)+t\to +\infty$ as $t\to +\infty$ and $s(t)+t\to-\infty$ as $t\to -\infty$. 
Together with \eqref{abcdeq} this yields 
$b=c=0$. Since $ad=1$ we can assume that $a>0, d>0$ and $a=e^{\tau/2},d=e^{-\tau/2}$ for some $\tau\in\R$. This implies that
$g^{-1}\gamma h=a_\tau$ or $y=\varphi_\tau(x)$. 
Finally, using $\big||a|-1\big|+\big||d|-1\big|<\eps_0$, we have
$e^{|\tau|/2}-e^{-|\tau|/2}<\eps_0=e^{\eps/2}-e^{-\eps/2}$, consequently $|\tau|< \eps$ which completes the proof. 
{\hfill $\Box$}

\subsection{Kinematic expansiveness, separation}
This subsection is devoted to demonstrate the  positive/negative strong kinematic expansiveness for the horocycle flow. 
It is also showed that the horocycle flow is not BW\,--\,expansive while the geodesic flow is not positive/negative 
separating. 
	\begin{definition}[\cite{artigue}, Kinematic expansive]\label{dnke}\rm Let $(M,d)$ be a compact metric space. 
		A continuous flow $\phi:\R\times M\longrightarrow M$ is called {\em kinematic expansive} if for each $\eps>0$, there exists $\delta>0$ with the following property. If
		\begin{equation}\label{exdn}d(\phi_t(x),\phi_{t}(y))<\delta \quad \mbox{for all}\quad t\in \R
		\end{equation}
		then $y=\phi_\tau(x)$ for some $\tau\in(-\eps,\eps)$.
	\end{definition}
If the inequality in  \eqref{exdn} holds for $t\in [0,\infty)$ (resp. $t\in (-\infty, 0]$) then the flow is called  `positive kinematic expansive' (resp. `negative kinematic expansive').
	If the condition $\tau\in(-\eps,\eps)$ is ignored, the flow is called  separating in the sense of Gura.
	\begin{definition}[\cite{gura}, Separating]\rm Let $(M,d)$ be a compact metric space. 
		A continuous flow $\phi:\R\times M\longrightarrow M$ is called {\em separating} if  there exists $\delta>0$ with the following property. If
		\begin{equation}\label{sdn}
		d(\phi_t(x),\phi_t(y))<\delta \quad \mbox{for all}\quad t\in \R
		\end{equation} 
		then $y=\phi_\tau(x)$ for some $\tau\in\R$; i.e. $x$ and $y$ lie on the same orbit. 
	\end{definition}
The number $\delta$ is called a `separating constant'.
If the inequality in  \eqref{sdn} holds only for $t\in [0,\infty)$ (resp. $t\in (-\infty, 0]$) then the flow is called  `positive separating' (resp. `negative separating').

\begin{definition}\label{tcdn}\rm Let $M$ be a metric space and $\phi,\psi:\R\times M\to M$ be continuous flows. We say 
	$(\phi_t)_{t\in\R}$ is a {\em time change} of $(\psi_t)_{t\in\R}$ if 
	for every $x\in M$ the orbits $\phi_\R(x)$, $\psi_\R(x)$ 
	and their orientations coincide.
\end{definition}

\begin{definition}\rm Let $M$ be a compact metric space. A continuous flow $\phi:\R\times M\longrightarrow M$ is called 
	{\em positive} (cores. {\em negative}) {\em strong separating} (cores. {\em strong expansive})
	if every time change of $(\phi_t)_{t\in\R}$ is {\em positive} (cores. {\em negative}) {\em separating} (cores. {\em expansive}).
\end{definition}
The following result is well-known:
\begin{proposition}[\cite{gura}]\label{ss} The horocycle flow on a compact surface of negative curvature is positive and negative strong separating.
\end{proposition}

The next result gives a stronger property of the horocycle flow  on $\Gamma\backslash\H^2$
which is a compact Riemann surface with constant negative curvature. 
	
	\begin{theorem}\label{pne}
		The horocycle flow $(\theta_t^\X)_{t\in\R}$ is positive and negative strong kinematic expansive.
	\end{theorem}
	{\bf Proof\,:}
	We consider the positive strong kinematic expansiveness only. Since the positive strong kinematic expansiveness is invariant under equivalence, 
	it follows from \eqref{cr} that it suffices to show that the flow $(\theta_t)_{t\in\R}$ is positive strong kinematic expansive.
	Let $(\psi_t)_{t\in\R}$ be a time change of the flow $(\theta_t)_{t\in\R}$ (see Definition \ref{tcdn}). 
	There exist  continuous functions $\alpha,\beta: \R\times X\to\R $ satisfying 
	\[\theta_t(x)=\psi_{\alpha(t,x)}(x)\quad\mbox{and}\quad \psi_t(x)=\theta_{\beta(t,x)}(x) \quad \mbox{for all}\quad (t, x)\in\R\times X. \]
	For every $\eps>0$, there is $\rho=\rho(\eps)\in (0,\eps_*)$ so that
	\begin{equation}\label{stt2} |\alpha(t,x)|<\eps\quad\mbox{for all}\quad (t,x)\in (-\rho,\rho)\times X;
	\end{equation} 
	recall $\eps_*$ from Lemma \ref{hy}. 
	With $\rho>0$ above, take $\delta=\delta(\rho)<\eps_*$ as in Lemma \ref{konvexa}\,(b). According to Proposition \ref{ss}, the flow $(\theta_t)_{t\in\R}$ is positive strong separating. This implies that the time change 
	$(\psi_t)_{t\in\R}$ is positive separating.
	Let $0<\varrho<\delta$ be a separating constant. For any $x,y\in X$, if 
	\begin{equation}\label{stt3}
	d_X(\psi_t(x),\psi_t(y))<\varrho\quad
	\mbox{for} \quad t\geq 0
	\end{equation}
	then $y=\psi_r(x)$ for some $r\in \R$. We need to show that $|r|<\eps$. There exists a unique $s\in \R$ such that $r=\alpha(s,x)$ and  $y=\theta_s(x)$.
	Due to \eqref{stt2}, it remains to verify that $|s|<\rho$.
	
	First, we rewrite \eqref{stt3} as 
	\[ d_X(\theta_{\beta(t,x)}(x),\theta_{\beta(t,y)+s}(x))<\varrho\quad\mbox{for}\quad t\geq 0. \]
	Denote $s_1(t)=\beta(t,x), s_2(t)=\beta(t,y)+s$ and write $x=\Gamma g$ for $g\in\PSL(2,\R)$. 
	For every $t\geq 0$, there is $\gamma(t)\in\Gamma$ so that
	\begin{equation*}\label{strong1}
	d_X(\theta_{s_1(t)}(x),\theta_{s_2(t)}(x))=d_\G(\gamma(t) gb_{s_1(t)}, gb_{s_2(t)})<\varrho. 
	\end{equation*} 
	Analogously to the argument in the proof of Theorem \ref{geoex}, we can check that $\gamma(t)=\gamma(0)=:\gamma$  for all $t\geq 0$ and then
	\[d_\G(\gamma gb_{s_1(t)}, gb_{s_2(t)})<\varrho\quad \mbox{for}\quad t\geq 0. \]
	or
	\[d_\G(b_{-s_2(t)}g^{-1}\gamma gb_{s_1(t)}, e)<\varrho\quad \mbox{for}\quad t\geq 0. \]
	According to Lemma \ref{konvexa}\,(b), there is $K=\big(\scriptsize\begin{array}{cc} k_{11}&k_{12}\\k_{21}&k_{22}\end{array}\big)\in\SL(2,\R)$ such that $g^{-1}\gamma g=\pi(K)$ and
	\begin{equation}\label{st1}|k_{11}-k_{21}s_2(t)-1 | + | (k_{11}-k_{21}s_2(t))s_1(t)-k_{22}s_2(t)+k_{12}| 
	+ | k_{21}|+ |k_{21}s_1(t)+k_{22}-1|<\rho, 
	\end{equation}
	using 
	\[ B_{-s_2(t)}KB_{s_1(t)}=\Big(\begin{array}{cc} k_{11}-k_{21}s_2(t)&(k_{11}-k_{21}s_2(t))s_1(t)-k_{22}s_2(t)+k_{12}\\
	k_{21}&k_{21}s_1(t)+k_{22}\end{array}\Big). \]
	Noting that $s_1(t)\to +\infty$ and $s_2(t)\to+\infty$ as $t\to +\infty$, 
	it follows from \eqref{st1}  that $k_{21}=0$ and
	\begin{equation}\label{st2} |k_{11}-1|+|k_{11}s_1(t)-k_{22}s_2(t)+k_{12}|+|k_{21}|+|k_{22}-1|<\rho. 
	\end{equation} 
	This yields 
	$\tr(\gamma)=\tr(g^{-1}\gamma g)=|k_{11}+k_{22}|<2+\rho<2+\eps_*$. 
	Owing to Lemma \ref{hy}, we have $\gamma=e$ and hence  $k_{11}=k_{22}=1, k_{12}=k_{21}=0$.
	Then \eqref{st2} implies $|s_1(t)-s_2(t)|<\rho$ for all $t\geq 0$. Finally, $s_1(0)=0, s_2(0)=s$ imply $|s|<\rho$ that completes the proof.
	{\hfill$\Box$}
		
\begin{remark}\label{rm}\rm 
	
	(a) The horocycle  flow is not BW\,--\,expansive. 
	Indeed, for any $\delta>0$, we need to find $x, y\in X$ and $s:\R\to \R$ continuous with $s(0)=0$ such that
	$d_X(\theta_t(x),\theta_{s(t)}(y))<\delta\quad\mbox{for all}\quad t\in\R$ but the orbits of $x$ and $y$ 
	do not coincide.  Take $\rho=\rho(\delta)$ as in Lemma \ref{konvexa}\,(a) and choose any $x=\Gamma g$ and $y=\Gamma h$ with $h,g\in\PSL(2,\R), h\ne g, h^{-1}g=\pi(K)$, $K=\Big(\begin{array}{cc}a&0\\ 0& d\end{array}\Big)$
	such that $ad=1$, $|a-1|<\rho, |d-1|<\rho$ and 
	$\tr(h^{-1}g)=|a+d|<2+\eps_*$; here $\eps_*>0$ is from Lemma \ref{hy}. We obtain
	$d_\G(h^{-1}g,e)<\delta$ due to Lemma \ref{konvexa}\,(a). Setting $s(t)=\frac{d}{a}t$, we have
	\begin{eqnarray*}
		d_X(\theta_{s(t)}(x),\theta_t(y))
		&=&d_X(\Gamma gb_{s(t)},\Gamma hb_t)\leq d_\G(gb_{s(t)},hb_t )=d_\G(b_{-t}h^{-1}gb_{s(t)},e)\\
		&=&d_\G(h^{-1}g,e)<\delta\quad\mbox{for all}\quad t\in\R;
	\end{eqnarray*}
	using $b_{-t}h^{-1}gb_{s(t)}=h^{-1}g$. It remains to verify that 
	$x$ and $y$ are not in the same orbit. Indeed, otherwise there would exist  $\tau\in\R$ such that $y=\theta_\tau(x)$, then $\gamma h= g b_\tau$ for some $\gamma\in\Gamma$ implies $\tr(\gamma)=\tr(gb_\tau h^{-1})=\tr(b_\tau h^{-1}g)=|a+d|<2+\eps_*$.
	It follows from Lemma \ref{hy} that $\gamma=e$. This yields 
	$b_{-\tau}=h^{-1}g=\pi(K)$ and hence $\tau=0, h=g$ which contradicts to $h\ne g$. 
	Therefore 	$x$ and $y$ do not lie in the same orbit and the horocycle flow is not BW\,--\,expansive. In fact, this flow is not geometric separating (see \cite[Definition 2.21]{artigue}).  
	
	\medskip
	(b) The geodesic flow is neither positive nor negative separating. Indeed, we consider the equivalent flow
	$(\varphi_t)_{t\in\R}$. Since the group $\Gamma$ is discrete,
	for every $\delta>0$, there is an $s\in (-\delta,\delta)$ such that
	 $a_tb_{-s}\notin\Gamma$ for all $t\in\R$. 
	Set $x=\Gamma e$ and $y=\Gamma b_{s}$ to have
	\begin{eqnarray*}
		d_X(\varphi_{t}(x),\varphi_t(y))
		=d_X(\Gamma a_{t},\Gamma b_{s}a_t)\leq d_\G(a_{t},b_sa_t )\leq  |s|e^{-t}	<\delta\quad \mbox{for all}\quad t\geq 0.
	\end{eqnarray*}
	However, if
	$y=\varphi_\tau(x)$ then $\Gamma b_s=\Gamma a_\tau$ implies that 
	there is $\gamma=a_\tau b_{-s}\in\Gamma$  which is a contradiction,
	whence $(\varphi_t)_{t\in\R}$
	is not positive separating. 
	In the same manner one obtains that the flow $(\varphi_t)_{t\in\R}$  is not negative separating.

	\medskip
	(c) It is worth mentioning that the geodesic flow is BW\,--\,expansive but neither positive nor negative kinematic expansive while the horocycle flows are positive and negative kinematic expansive but not BW\,--\,expansive.
	{\hfill $\diamondsuit$}

\end{remark}

\subsection{KH\,--\, expansiveness}
In \cite{kh} Katok and Hasselblatt introduce the following expansiveness:
\begin{definition}[\cite{kh}, {KH}\,--\,expansive]\rm
	Let $(M,d)$ be a compact space. A continuous flow $\phi_t:M\longrightarrow M$ is called {\em KH\,--\,expansive} if there exists $\delta>0$ with the following property. If $x\in X, s:\R\to\R$ is continuous, $s(0)=0$ and $d(\varphi_t(x),\varphi_{s(t)}(x))<\delta
	$ for all $t\in\R, y\in X$ is such that $d(\varphi_t(x),\varphi_{s(t)}(y))<\delta$ for all $t\in\R$
	then $x$ and $y$ lie on the same orbit.
\end{definition}
It is clear that {KH}\,--\,expansiveness is weaker than BW\,--\,expansiveness but implies separation. 
Furthermore, one has the following result:	
\begin{proposition}[\cite{artigue2}]\label{KHthm} A flow  on a compact metric space is KH\,--\,expansive if and only if it is separating and the set of its fixed points is open.
\end{proposition}

It follows immediately from propositions \ref{fixedpoint} and \ref{KHthm}  that the  flow $(\theta_t)_{t\in\R}$ is KH\,--\,expansive, and hence the horocycle flow $(\theta_t^\X)_{t\in\R}$ is KH\,--\,expansive owing to \eqref{cr}.
Nevertheless we can verify it directly.
\begin{theorem}\label{khthm}
	The horocycle flow $(\theta^\X_t)_{t\in\R}$ is KH\,--\,expansive.
\end{theorem}	
\noindent 
{\bf Proof\,:}   If $x,y\in X, s:\R\to\R$ is continuous, $s(0)=0$ and 
\begin{equation}\label{s1}
d_X(\theta_t(x),\theta_{s(t)}(x))<\delta
\quad\mbox{for all}\quad  t\in\R
\end{equation}
and
\begin{equation}\label{s2}
d_X(\theta_t(x),\theta_{s(t)}(y))<\delta\quad\mbox{for all}\quad t\in\R.
\end{equation} 
Analogously to the proof of Theorem \ref{geoex}, using \eqref{s1} we can show that 
there is  $M>0$ such that 
\[|s(t)-t|<M\quad\mbox{for all}\quad t\in\R. \]
This means that \begin{equation}\label{s3}
s(t)\to +\infty\quad \mbox{as}\quad t\to +\infty.
\end{equation} 
It follows from \eqref{s1} and \eqref{s2} that
\[d_X(\theta_{s(t)}(x), \theta_{s(t)}(y))<2\delta\quad\mbox{for all}\quad t\in\R.\]
Together with \eqref{s3},  this follows in the same manner of the proof of Theorem \ref{pne}.
{\hfill$\Box$}

\end{document}